\documentclass[11pt,letterpaper]{article}
\usepackage[T1]{fontenc}
\usepackage[utf8]{inputenc}
\usepackage{lmodern}
\usepackage{amsmath,amssymb,mathtools}
\usepackage{graphicx}
\usepackage{booktabs,array,tabularx}
\usepackage[margin=1in]{geometry}
\usepackage[labelfont=bf,font=small]{caption}
\usepackage{microtype}
\usepackage{enumitem}
\usepackage{fancyhdr}
\usepackage[colorlinks=true,linkcolor=black,urlcolor=blue!50!black,citecolor=black]{hyperref}
\usepackage{xcolor}
\usepackage{float}
\renewcommand{\arraystretch}{1.14}
\setlist[enumerate]{leftmargin=*,itemsep=0.4em,topsep=0.3em}
\title{\textbf{Alternative Tikhonov-Type Regularization for\\the Calder\'on Inverse Conductivity Problem}}
\author{Pavel B. Dubovski, Charles Morehead}
\date{}
\begin{document}
\maketitle
\vspace{-1.1em}
\begin{abstract}

This paper investigates the application of two regularization methods, classical Tikhonov and Alternative Tikhonov regularization, to the Calderón inverse conductivity problem. We develop an iterative reconstruction scheme based on linearizing the Dirichlet-to-Neumann map and approximating conductivity updates using Gaussian basis functions. The resulting linear systems are solved using both regularization methods, with particular attention given to the derivation and numerical implementation of Alternative Tikhonov regularization. We compare the performance of these methods through numerical experiments on a synthetic conductivity containing two localized regions of differing conductivity. A total of 666 reconstruction tests are performed across varying regularization parameters, damping factors, and iteration counts. The lowest relative reconstruction errors obtained are 0.064730 for Alternative Tikhonov and 0.065173 for classical Tikhonov regularization. Although both methods produce comparable reconstructions, Alternative Tikhonov achieves the lowest observed error. The results demonstrate the effectiveness of both approaches within the proposed reconstruction scheme and motivate further investigation into parameter selection, numerical stability, and reconstruction accuracy. 
\end{abstract}
\section{Introduction}
The Calder\'on inverse problem, also known as electrical impedance tomography in the medical field, seeks to recover an unknown conductivity function $\gamma:\Omega\to \mathbb{R}^+$ with bounded domain $\Omega\subset \mathbb{R}^n$ from certain boundary measurements. With $u:\Omega\to \mathbb{R}$ representing the electric potential, the forward problem is:
\begin{equation*}
\nabla\cdot(\gamma\nabla u)=0\quad\text{in }\Omega.
\end{equation*}
\begin{equation*}
u|_{\partial\Omega}=f.
\end{equation*}
We are also given the output of the Dirichlet-to-Neumann map
\begin{equation*}
\Lambda_\gamma:f\mapsto\left.\gamma\frac{\partial u}{\partial n}\right|_{\partial\Omega}.
\end{equation*}
We are given this information to find $\gamma$ \cite{Calderon}.

 This is an ill-posed problem. Small perturbations in the Dirichlet-to-Neumann map can lead to large differences in the reconstruction of $\gamma$ \cite{Uhlmann}. Since the problem is also non-linear, a stable reconstruction scheme must provide a linearization of the problem and a regularization method to solve that linear system. In order to obtain a linearized system, we perform multiple experiments with different boundary measurements $f_i$, resulting in multiple potentials $u_i$. To obtain stable approximations of the conductivity, we develop an iterative reconstruction algorithm incorporating two regularization methods: classical Tikhonov regularization and Alternative Tikhonov (AT) regularization. We then compare their performance in reconstructing the unknown conductivity.

\section{Iterative Reconstruction Scheme}
We construct an iterative approximation to the conductivity $\gamma$ as follows.
\begin{enumerate}
\item Choose an initial approximation $\gamma_0$.
\item Given $\gamma^{(k)}$, solve the forward problem
\begin{equation*}
\nabla\cdot(\gamma^{(k)}\nabla u^{(k)})=0\quad\text{in }\Omega,\qquad u^{(k)}|_{\partial\Omega}=f.
\end{equation*}
\item Solve the adjoint problem
\begin{equation*}
-\nabla\cdot(\gamma^{(k)}\nabla v^{(k)})=0,\qquad v^{(k)}|_{\partial\Omega}=g.
\end{equation*}
\item Find the deviation $\delta\gamma^{(k)}$ of $\gamma^{(k)}$:
\begin{equation*}
\gamma^{(k+1)}=\gamma^{(k)}+\delta\gamma^{(k)}.
\end{equation*}
\end{enumerate}
The Dirichlet-to-Neumann map is nonlinear in $\gamma$, so we linearize about $\gamma^{(k)}$ \cite{Hyvonen}:
\begin{equation*}
\Lambda_{\gamma^{(k)}+\delta\gamma^{(k)}}(f)
=
\Lambda_{\gamma^{(k)}}(f)
+
D\Lambda_{\gamma^{(k)}}[\delta\gamma^{(k)}](f)
+
\text{higher-order terms}.
\end{equation*}

Let
\begin{equation*}
\gamma_\varepsilon^{(k)}
=
\gamma^{(k)}+\varepsilon\delta\gamma^{(k)},
\qquad
u_\varepsilon
=
u^{(k)}+\varepsilon w+o(\varepsilon),
\end{equation*}
where
\begin{equation*}
\nabla\cdot(\gamma_\varepsilon^{(k)}\nabla u_\varepsilon)=0,
\qquad
u_\varepsilon|_{\partial\Omega}=f.
\end{equation*}

Substituting these expressions into the forward equation gives
\begin{equation*}
\nabla\cdot\left[
(\gamma^{(k)}+\varepsilon\delta\gamma^{(k)})
\nabla(u^{(k)}+\varepsilon w+o(\varepsilon))
\right]=0.
\end{equation*}

Expanding, we obtain
\begin{equation*}
\begin{aligned}
\nabla\cdot(\gamma^{(k)}\nabla u^{(k)}) + \varepsilon\nabla\cdot(\gamma^{(k)}\nabla w) + \varepsilon\nabla\cdot(\delta\gamma^{(k)}\nabla u^{(k)}) + o(\varepsilon)=0.
\end{aligned}
\end{equation*}

Since
\begin{equation*}
\nabla\cdot(\gamma^{(k)}\nabla u^{(k)})=0,
\end{equation*}
keeping only the first-order terms gives
\begin{equation*}
\nabla\cdot(\gamma^{(k)}\nabla w)
=
-\nabla\cdot(\delta\gamma^{(k)}\nabla u^{(k)}),
\qquad
w|_{\partial\Omega}=0.
\end{equation*}

We now consider the Dirichlet-to-Neumann map
\begin{equation*}
\Lambda_\gamma(f)
=
\left.\gamma\frac{\partial u}{\partial n}
\right|_{\partial\Omega}.
\end{equation*}

Substituting the perturbed conductivity and potential gives
\begin{equation*}
\begin{aligned}
\Lambda_{\gamma_\varepsilon^{(k)}}(f)
={}&
\left.
(\gamma^{(k)}+\varepsilon\delta\gamma^{(k)})\right.\\
&\left.
\cdot\left(
\frac{\partial u^{(k)}}{\partial n}
+\varepsilon\frac{\partial w}{\partial n}
\right)\right|_{\partial\Omega}
+o(\varepsilon).
\end{aligned}
\end{equation*}

Expanding and keeping first-order terms, we have
\begin{equation*}
\begin{aligned}
\Lambda_{\gamma_\varepsilon^{(k)}}(f)
={}&
\Lambda_{\gamma^{(k)}}(f)\\
&+\varepsilon\left(
\delta\gamma^{(k)}
\frac{\partial u^{(k)}}{\partial n}
+
\gamma^{(k)}
\frac{\partial w}{\partial n}
\right)
+o(\varepsilon).
\end{aligned}
\end{equation*}

Therefore,
\begin{equation*}
D\Lambda_{\gamma^{(k)}}[\delta\gamma^{(k)}](f)
=
\delta\gamma^{(k)}
\frac{\partial u^{(k)}}{\partial n}
+
\gamma^{(k)}
\frac{\partial w}{\partial n}.
\end{equation*}

If $v^{(k)}$ is the adjoint state, then
\begin{equation*}
\begin{aligned}
\int_{\partial\Omega}
v^{(k)}
D\Lambda_{\gamma^{(k)}}[\delta\gamma^{(k)}](f)
\,dS
=
\int_{\partial\Omega}
v^{(k)}\delta\gamma^{(k)}
\frac{\partial u^{(k)}}{\partial n}\,dS
+
\int_{\partial\Omega}
v^{(k)}\gamma^{(k)}
\frac{\partial w}{\partial n}\,dS.
\end{aligned}
\end{equation*}

We now use integration by parts to simplify the second boundary integral. By the divergence theorem,
\begin{equation*}
\begin{aligned}
\int_{\partial\Omega}
v^{(k)}\gamma^{(k)}
\frac{\partial w}{\partial n}\,dS =
\int_\Omega
\gamma^{(k)}
\nabla v^{(k)}\cdot\nabla w\,dx +
\int_\Omega
v^{(k)}
\nabla\cdot(\gamma^{(k)}\nabla w)\,dx.
\end{aligned}
\end{equation*}

Since $v^{(k)}$ satisfies the homogeneous adjoint equation
\begin{equation*}
\nabla\cdot(\gamma^{(k)}\nabla v^{(k)})=0,
\end{equation*}
and $w|_{\partial\Omega}=0$, integration by parts gives
\begin{equation*}
\begin{aligned}
\int_\Omega
\gamma^{(k)}
\nabla v^{(k)}\cdot\nabla w\,dx =
\int_{\partial\Omega}
w\gamma^{(k)}
\frac{\partial v^{(k)}}{\partial n}\,dS - \int_\Omega w\nabla\cdot(\gamma^{(k)}\nabla v^{(k)})\,dx =0.
\end{aligned}
\end{equation*}

Thus,
\begin{equation*}
\begin{aligned}
\int_{\partial\Omega}
v^{(k)}\gamma^{(k)}
\frac{\partial w}{\partial n}\,dS = \int_\Omega v^{(k)} \nabla\cdot(\gamma^{(k)}\nabla w)\,dx.
\end{aligned}
\end{equation*}

Using the linearized equation for $w$, we obtain
\begin{equation*}
\begin{aligned}
\int_{\partial\Omega}
v^{(k)}\gamma^{(k)}
\frac{\partial w}{\partial n}\,dS = -\int_\Omega v^{(k)}\nabla\cdot(\delta\gamma^{(k)}\nabla u^{(k)})\,dx.
\end{aligned}
\end{equation*}

Integrating by parts once more,
\begin{equation*}
\begin{aligned}
\int_{\partial\Omega}
v^{(k)}\gamma^{(k)}
\frac{\partial w}{\partial n}\,dS=
-\int_{\partial\Omega}
v^{(k)}\delta\gamma^{(k)}
\frac{\partial u^{(k)}}{\partial n}\,dS+\int_\Omega\delta\gamma^{(k)}\nabla u^{(k)}\cdot\nabla v^{(k)}\,dx.
\end{aligned}
\end{equation*}

Substituting this into our expression for the derivative of the Dirichlet-to-Neumann map, the two boundary terms cancel. Therefore,
\begin{equation*}
\begin{aligned}
\int_{\partial\Omega}
v^{(k)}
D\Lambda_{\gamma^{(k)}}[\delta\gamma^{(k)}](f)
\,dS =
\int_\Omega
\delta\gamma^{(k)}
\nabla u^{(k)}\cdot\nabla v^{(k)}\,dx.
\end{aligned}
\end{equation*}

We follow \cite{Hyvonen} and approximate deviation of $\Lambda$
\begin{equation*}
D\Lambda_{\gamma^{(k)}}[\delta\gamma](f)=\Lambda_{\mathrm{meas}}(f)-\Lambda_{\gamma^{(k)}}(f).
\end{equation*}
Thus we obtain:
\begin{equation}\label{eq:one}
\int_\Omega\delta\gamma^{(k)}\nabla u^{(k)}\cdot\nabla v^{(k)}\,dx\approx\int_{\partial\Omega}v^{(k)}\left(\Lambda_{\mathrm{meas}}(f)-\Lambda_{\gamma^{(k)}}(f)\right)\,dS.
\end{equation}
In order to solve the relationship (1) for $\delta\gamma^{(k)}$, we discretize (1). We approximate a $\delta\gamma$ by substituting
\begin{equation*}
\delta\gamma^{(k)}=\sum_{j=1}^M c_j\psi_j,
\end{equation*}
where basis functions $\psi_j$ are Gaussian functions
\begin{equation*}
\psi_j(x,y)=\exp\left(-\frac{(x-x_j)^2+(y-y_j)^2}{2\sigma^2}\right)
\end{equation*}
where constant $\sigma$ and centers $(x_j,y_j)$ in a grid formation are chosen such that the whole domain $\Omega$ is covered.
Then using the relationship in (1) we have a linear system to solve for coefficients $c_j$. Let us call this system
\begin{equation}\label{eq:two}
Ac=r.
\end{equation}
Suppose we choose $M$ basis functions $\psi_j$ and use $N$ pairs of experimental boundary inputs $f_i$ and control boundary inputs $g_i$ from which we can find potentials $u_i^{(k)}$ and $v_i^{(k)}$ respectively. Thus building the following system:
\begin{equation*}
A\in\mathbb R^{N\times M},\quad c\in\mathbb R^M,\quad r\in\mathbb R^N
\end{equation*}
\begin{equation*}
A_{i,j}=\int_\Omega\psi_j\nabla u_i^{(k)}\cdot\nabla v_i^{(k)}\,dx,\qquad r_i=\int_{\partial\Omega}v_i^{(k)}\left(\Lambda_{\mathrm{meas}}(f_i)-\Lambda_{\gamma^{(k)}}(f_i)\right)\,dS.
\end{equation*}
The linear system (2) is ill-posed, so we solve for $c$ using regularization techniques. The regularization techniques used will be described in section 3.
\begin{enumerate}[start=5]
\item Compute relative error
\begin{equation*}
\frac{\|\gamma^{(k+1)}-\gamma_{\mathrm{true}}\|}{\|\gamma_{\mathrm{true}}\|}
\end{equation*}
Experience shows that after sufficiently many $n$ iterations error will begin to increase due to overfitting.
\item Repeat from step 2 or stop iteration if error is observed to begin to increase.
\end{enumerate}

\section{Regularization}
We will investigate two methods of regularization. The first of which is classical Tikhonov regularization, in which we find
\begin{equation*}
\min_{c\in\mathbb R^M}\|Ac-r\|^2+\alpha\|c\|^2.
\end{equation*}
Where $\alpha>0$ is a constant parameter. The second of which is a lesser known regularization method established by A.N. Tikhonov \cite{Tikhonov} and further developed by E.E. Tyrtyshnikov \cite{Tyrtyshnikov}. This method involves solving for constant $\rho$ in the fixed point problem
\begin{equation*}
\rho=\delta+\mu g(\rho)
\end{equation*}
where $\delta,\mu>0$ are constant parameters and
\begin{equation*}
g(\rho)=\min_{c\in\mathbb R^M}\|c\|\quad\text{subject to }\|Ac-r\|=\rho.
\end{equation*}
The solution is the vector $c$ satisfying the above minimum. In parallel with testing the quality of the results from reconstruction method, we will also compare and evaluate the effectiveness of each regularization method. We will refer to this as Alternative Tikhonov (AT) regularization, whereas we will refer to the first method as classical Tikhonov regularization.

\subsection{Classical Tikhonov Regularization}
We now describe the classical Tikhonov regularization method for solving the linear system
\begin{equation*}
Ac=r.
\end{equation*}
We solve the classical Tikhonov regularization problem \cite{Tikhonov}
\begin{equation*}
\min_{c\in\mathbb R^M}\bigl(\|Ac-r\|^2+\alpha\|c\|^2\bigr),
\end{equation*}
where $\alpha>0$ is the regularization parameter.
The parameter $\alpha$ controls the strength of the regularization. If $\alpha$ is small, the method behaves more like ordinary least squares and tries harder to fit the projected system $Ac=r$. If $\alpha$ is large, the method penalizes large coefficient vectors more strongly, producing a more stable but possibly more smoothed update. We will observe the effects of $\alpha$ in the results and see which value is best suited to this problem.

\noindent\textbf{Algorithm}

This will serve as a look into how step 4 from section 2 is performed using classical Tikhonov regularization.
\begin{enumerate}
\item Form linear system $A$ and $r$ from given data.
\item Solve for the regularized coefficient vector
\begin{equation*}
c=(A^TA+\alpha I)^{-1}A^Tr.
\end{equation*}
\item Reconstruct the conductivity update by
\begin{equation*}
\delta\gamma=\sum_{j=1}^M(c)_j\psi_j.
\end{equation*}
\item Update the conductivity by
\begin{equation*}
\gamma^{(k+1)}=\max\{\gamma^{(k)}+\omega\delta\gamma^{(k)},\gamma_{\min}\}.
\end{equation*}
In order to avoid too large steps we introduce damping parameter $\omega>0$ and $\gamma_{\min}>0$ ensures that $\gamma^{(k+1)}$ remains positive.
\end{enumerate}

\subsection{Alternative Tikhonov (AT) regularization}
In conducting Alternative Tikhonov regularization our main focus is solving the fixed point problem [8, 9]
\begin{equation*}
\rho=\delta+\mu g(\rho)
\end{equation*}
Consider the equivalent constrained problem
\begin{equation*}
g^2(\rho)=\min\|c\|^2\quad\text{subject to}\quad\|Ac-r\|^2=\rho^2.
\end{equation*}
We form the Lagrangian
\begin{equation*}
L(c,\lambda)=\|c\|^2+\lambda\bigl(\|Ac-r\|^2-\rho^2\bigr).
\end{equation*}
Taking the derivative with respect to $c$, we obtain
\begin{equation*}
\nabla_c L=2c+2\lambda A^T(Ac-r)=0.
\end{equation*}
Therefore,
\begin{equation*}
c+\lambda A^TAc=\lambda A^Tr,
\end{equation*}
so
\begin{equation*}
(I+\lambda A^TA)c=\lambda A^Tr.
\end{equation*}
Thus, for each $\lambda>0$, define
\begin{equation*}
c_\lambda=\lambda(I+\lambda A^TA)^{-1}A^Tr.
\end{equation*}
The Alternative Tikhonov condition requires
\begin{equation*}
\|Ac_\lambda-r\|=\delta+\mu\|c_\lambda\|.
\end{equation*}
Therefore, instead of solving for $\rho$, we solve the root problem
\begin{equation*}
H(\lambda)=0,
\end{equation*}
for $\lambda$ using Newton's method where
\begin{equation*}
H(\lambda)=\|Ac_\lambda-r\|-\delta-\mu\|c_\lambda\|.
\end{equation*}

\noindent\textbf{Algorithm}

This is how step 4 from section 2 is performed using Alternative Tikhonov regularization.
Given $A$, $r$, $\mu$, and $\delta$:
\begin{enumerate}
\item Define
\begin{equation*}
c_\lambda=\lambda(I+\lambda A^TA)^{-1}A^Tr.
\end{equation*}
\item Define the scalar function
\begin{equation*}
H(\lambda)=\|Ac_\lambda-r\|-\delta-\mu\|c_\lambda\|.
\end{equation*}
\item Find $\lambda^*>0$ such that
\begin{equation}\label{eq:three}
H(\lambda^*)=0.
\end{equation}
Since $I+\lambda A^TA$ is invertible for every $\lambda>0$, the vector $c_\lambda$ is well-defined for $\lambda>0$.
To solve the scalar equation
\begin{equation*}
H(\lambda)=0,
\end{equation*}
we use Newton's method. First define
\begin{equation*}
R(\lambda)=\|Ac_\lambda-r\|,\qquad C(\lambda)=\|c_\lambda\|.
\end{equation*}
Then
\begin{equation*}
H(\lambda)=R(\lambda)-\delta-\mu C(\lambda).
\end{equation*}
Differentiating gives
\begin{equation*}
H'(\lambda)=R'(\lambda)-\mu C'(\lambda).
\end{equation*}
Since
\begin{equation*}
R(\lambda)=\|Ac_\lambda-r\|,
\end{equation*}
we have
\begin{equation*}
R'(\lambda)=\frac{(Ac_\lambda-r)^TAc'_\lambda}{\|Ac_\lambda-r\|}.
\end{equation*}
Similarly, since
\begin{equation*}
C(\lambda)=\|c_\lambda\|,
\end{equation*}
we have
\begin{equation*}
C'(\lambda)=\frac{c_\lambda^Tc'_\lambda}{\|c_\lambda\|}.
\end{equation*}
Therefore
\begin{equation*}
H'(\lambda)=\frac{(Ac_\lambda-r)^TAc'_\lambda}{\|Ac_\lambda-r\|}-\mu\frac{c_\lambda^Tc'_\lambda}{\|c_\lambda\|}.
\end{equation*}
It remains to justify that $H(\lambda)$ is decreasing. Using the singular value decomposition
\begin{equation*}
A=U\Sigma V^T,\qquad\beta=U^Tr,
\end{equation*}
the regularized solution can be written as
\begin{equation*}
c_\lambda=\sum_i\frac{\lambda\sigma_i}{1+\lambda\sigma_i^2}\beta_i v_i.
\end{equation*}
Thus
\begin{equation*}
R(\lambda)^2=\|Ac_\lambda-r\|^2=\sum_i\left(\frac{\beta_i}{1+\lambda\sigma_i^2}\right)^2+\|r_\perp\|^2,
\end{equation*}
where $r_\perp$ is the component of $r$ orthogonal to the range of $A$. Differentiating gives
\begin{equation*}
\frac{d}{d\lambda}R(\lambda)^2=-2\sum_i\frac{\beta_i^2\sigma_i^2}{(1+\lambda\sigma_i^2)^3}\le0.
\end{equation*}
Since
\begin{equation*}
\frac{d}{d\lambda}R(\lambda)^2=2R(\lambda)R'(\lambda),
\end{equation*}
it follows that
\begin{equation*}
R'(\lambda)\le0.
\end{equation*}
Similarly,
\begin{equation*}
C(\lambda)^2=\|c_\lambda\|^2=\sum_i\left(\frac{\lambda\sigma_i}{1+\lambda\sigma_i^2}\beta_i\right)^2.
\end{equation*}
Therefore
\begin{equation*}
C(\lambda)^2=\sum_i\beta_i^2\frac{\lambda^2\sigma_i^2}{(1+\lambda\sigma_i^2)^2}.
\end{equation*}
Differentiating gives
\begin{equation*}
\frac{d}{d\lambda}C(\lambda)^2=2\sum_i\frac{\lambda\beta_i^2\sigma_i^2}{(1+\lambda\sigma_i^2)^3}\ge0
\end{equation*}
for $\lambda>0$. Since
\begin{equation*}
\frac{d}{d\lambda}C(\lambda)^2=2C(\lambda)C'(\lambda),
\end{equation*}
we obtain
\begin{equation*}
C'(\lambda)\ge0.
\end{equation*}
Thus
\begin{equation*}
H'(\lambda)=R'(\lambda)-\mu C'(\lambda)\le0.
\end{equation*}
 
Therefore $H(\lambda)$ is non-increasing. Furthermore,
if $A^Tr\neq 0$, then at least one term in
\begin{equation*}
\frac{d}{d\lambda}R(\lambda)^2
=
-2\sum_i
\frac{\beta_i^2\sigma_i^2}
{(1+\lambda\sigma_i^2)^3}
\end{equation*}
is strictly negative. Thus $R'(\lambda)<0$,
and since $C'(\lambda)\geq 0$, we obtain
\begin{equation*}
H'(\lambda)=R'(\lambda)-\mu C'(\lambda)<0.
\end{equation*}
Therefore $H(\lambda)$ is strictly decreasing
when $A^Tr\neq 0$. If $A^Tr=0$, we obtain trivial solution $c_\lambda=0$, so we will assume $A^Tr\neq 0$.

The Newton iteration is therefore
\begin{equation*}
\lambda_{j+1}=\lambda_j-\frac{H(\lambda_j)}{H'(\lambda_j)}.
\end{equation*}
The existence result for the approximate system guarantees a solution $\lambda$ to $H(\lambda)=0$ \cite{Tikhonov}\cite{Tyrtyshnikov}.

Since $H$ is strictly decreasing, any root of $H(\lambda)=0$ is unique. We now show that this root must be positive.

Recall that
\begin{equation*}
H(\lambda)=\|Ac_\lambda-r\|-\delta-\mu\|c_\lambda\|,
\end{equation*}
where
\begin{equation*}
c_\lambda=\lambda(I+\lambda A^TA)^{-1}A^Tr.
\end{equation*}

At $\lambda=0$, we have $c_0=0$. Therefore,
\begin{equation*}
\begin{aligned}
H(0)
&=\|Ac_0-r\|-\delta-\mu\|c_0\|\\
&=\|r\|-\delta.
\end{aligned}
\end{equation*}

Assuming $\delta<\|r\|$, we obtain
\begin{equation*}
H(0)>0.
\end{equation*}

Since $H$ is decreasing for $\lambda\geq 0$, and $H(0)>0$, any root $\lambda^*$ satisfying
\begin{equation*}
H(\lambda^*)=0
\end{equation*}
must satisfy
\begin{equation*}
\lambda^*>0.
\end{equation*}

Thus, assuming a root exists, it is unique and strictly positive.

We can therefore have the lower bracket as
\begin{equation*}
\lambda_{\mathrm{low}}=0.
\end{equation*}
For the upper bracket, we start with
\begin{equation*}
\lambda_{\mathrm{high}}=1
\end{equation*}
and repeatedly double it until
\begin{equation*}
H(\lambda_{\mathrm{high}})<0.
\end{equation*}
Once the interval
\begin{equation*}
[\lambda_{\mathrm{low}},\lambda_{\mathrm{high}}]
\end{equation*}
brackets the root, Newton's method can be applied safely.
\item Set
\begin{equation*}
c^*=c_{\lambda^*}.
\end{equation*}

\item Reconstruct the conductivity update by
\begin{equation*}
\delta\gamma=\sum_{j=1}^M c_j^*\psi_j.
\end{equation*}
\item Update the conductivity by
\begin{equation*}
\gamma^{(k+1)}=\max\{\gamma^{(k)}+\omega\delta\gamma^{(k)},\gamma_{\min}\}.
\end{equation*}
\end{enumerate}

\section{Numerical Results}
To study the effectiveness of our methods we test the inverse conductivity method on the square domain;
\begin{equation*}
\Omega=[0,1]^2
\end{equation*}
The unknown conductivity is a function
\begin{equation*}
\gamma(x,y),\;\; x,y\in \Omega
\end{equation*}
and the forward problem is
\begin{equation*}
-\nabla\cdot(\gamma\nabla u)=0\quad\text{in }\Omega,\qquad u|_{\partial\Omega}=f.\tag{$*$}
\end{equation*}
For different chosen boundary input $f$, we assume a benchmark conductivity $\gamma$ and solve $(*)$ for $u$. We then use $u$ to find boundary flux:
\begin{equation*}
q=\gamma\frac{\partial u}{\partial n}.
\end{equation*}
This boundary flux plays the role of the measured output data. We then take boundary flux $q$ and boundary potential $f$ as known and solve the inverse problem for $\gamma$, to be compared to the benchmark $\gamma$ for accuracy.
The code first chooses a true conductivity $\gamma_{\mathrm{true}}$ and generates boundary flux from it. As a benchmark we consider:
\begin{equation*}
\gamma_{\mathrm{true}}(x,y)=1+0.9e^{-\frac{(x-x_1)^2+(y-y_1)^2}{\sigma_1}}-0.4e^{-\frac{(x-x_2)^2+(y-y_2)^2}{\sigma_2}},
\end{equation*}
where
\begin{equation*}
x_1=0.68,\ y_1=0.36,\ x_2=0.32,\ y_2=0.72
\end{equation*}
and
\begin{equation*}
\sigma_1=0.012,\ \sigma_2=0.020
\end{equation*}
with a lower cutoff $\gamma_{\min}=0.2$ to keep the conductivity positive. Thus the true conductivity consists of a constant background near 1, one localized region of higher conductivity, and one localized region of lower conductivity.
The reconstruction is initialized with the constant guess
\begin{equation*}
\gamma_0(x,y)=1.
\end{equation*}
Starting from this guess we follow the algorithm as described in section 2 to try and reconstruct the true conductivity as accurately as possible. For each reconstruction, we run eight numerical experiments $u_i|_{\partial\Omega}=f_i$. From the adjoint $-\nabla\cdot(\gamma\nabla u)=0$ we choose for simplicity and six boundary functions $v_i|_{\partial\Omega}=g_i$. Each forward boundary function $f_i$ is paired with each test boundary function $g_i$ exactly once, giving $N=48$ boundary-function pairs at each experiment. The forward boundary functions are
\begin{align*}
f_1&=\begin{cases}\sin(\pi y),&\text{on the left boundary},\\0,&\text{otherwise},\end{cases}
&f_2&=\begin{cases}\sin(\pi y),&\text{on the right boundary},\\0,&\text{otherwise},\end{cases}\\[1em]
f_3&=\begin{cases}\sin(\pi x),&\text{on the bottom boundary},\\0,&\text{otherwise},\end{cases}
&f_4&=\begin{cases}\sin(\pi x),&\text{on the top boundary},\\0,&\text{otherwise},\end{cases}\\[1em]
f_5&=\begin{cases}\sin(2\pi y),&\text{on the left boundary},\\0,&\text{otherwise},\end{cases}
&f_6&=\begin{cases}\sin(2\pi x),&\text{on the bottom boundary},\\0,&\text{otherwise}.\end{cases}
\end{align*}
\begin{align*}
f_7&=\begin{cases}y(1-y),&\text{on the left boundary},\\-y(1-y),&\text{on the right boundary},\\0,&\text{otherwise},\end{cases}\\[1em]
f_8&=\begin{cases}x(1-x),&\text{on the bottom boundary},\\-x(1-x),&\text{on the top boundary},\\0,&\text{otherwise}.\end{cases}
\end{align*}
The six test boundary functions are
\begin{align*}
g_1&=\begin{cases}\sin(\pi y),&\text{on the left boundary},\\0,&\text{otherwise},\end{cases}
&g_2&=\begin{cases}\sin(\pi y),&\text{on the right boundary},\\0,&\text{otherwise},\end{cases}\\[1em]
g_3&=\begin{cases}\sin(\pi x),&\text{on the bottom boundary},\\0,&\text{otherwise},\end{cases}
&g_4&=\begin{cases}\sin(\pi x),&\text{on the top boundary},\\0,&\text{otherwise},\end{cases}
\end{align*}
\begin{align*}
g_5&=\begin{cases}y(1-y),&\text{on the left boundary},\\-y(1-y),&\text{on the right boundary},\\0,&\text{otherwise},\end{cases}\\[1em]
g_6&=\begin{cases}x(1-x),&\text{on the bottom boundary},\\-x(1-x),&\text{on the top boundary},\\0,&\text{otherwise}.\end{cases}
\end{align*}
 These same boundary functions are used for every reconstruction so that the effects of the regularization method and the tested parameters can be compared directly.
We will observe how the reconstruction reacts to the following parameters:
\begin{enumerate}
\item $n$: number of iterations,
\item $\omega$: damping factor for each update,
\item $\alpha$: constant for classical Tikhonov regularization,
\item $\delta,\mu$: constants for Alternative Tikhonov regularization.
\end{enumerate}
The values tested are
\begin{equation*}
n\in\{3,4,5,10,20,30\},
\end{equation*}
and
\begin{equation*}
\omega\in\{0.5,0.8,1.0\}.
\end{equation*}
For classical Tikhonov regularization, we also vary the regularization parameter $\alpha$, using
\begin{equation*}
\alpha\in\{10^{-5},3\cdot10^{-5},10^{-4},3\cdot10^{-4},10^{-3},3\cdot10^{-3},10^{-2}\}.
\end{equation*}
Therefore, the number of classical Tikhonov tests we did is
\begin{equation*}
6\cdot3\cdot7=126.
\end{equation*}
For Alternative Tikhonov regularization, we vary the parameters $\mu$ and $\delta$, using
\begin{equation*}
\mu\in\{10^{-3},3\cdot10^{-3},10^{-2},3\cdot10^{-2},10^{-1},3\cdot10^{-1}\},
\end{equation*}
and
\begin{equation*}
\delta\in\{10^{-4},3\cdot10^{-4},10^{-3},3\cdot10^{-3},10^{-2}\}.
\end{equation*}
Therefore, the number of Alternative Tikhonov tests we did is
\begin{equation*}
6\cdot3\cdot6\cdot5=540.
\end{equation*}
Altogether, the total number of reconstruction tests is
\begin{equation*}
126+540=666.
\end{equation*}
More detailed numerical results are given in the addendum. The Python implementation used for the numerical experiments is available on GitHub\cite{Github}. The following image shows the best result from each regularization method.

\section{Conclusion}
We developed a novel method of solving the inverse Calder\'on problem arising in medical tomography. To increase reliability, we applied classical Tikhonov regularization and compared with Alternative Tikhonov regularization inside the linearized Calder\'on reconstruction scheme. As a byproduct, we showed that, in this case, the Alternative Tikhonov scheme produced better results than classical Tikhonov regularization. Both methods were tested on the same benchmark conductivity, using the same boundary data, basis functions, initial guess, and update rule.
The best overall result came from Alternative Tikhonov regularization with
\begin{equation*}
n=3,\qquad\omega=0.8,\qquad\mu=10^{-2},\qquad\delta=10^{-4},
\end{equation*}
which gives relative reconstruction error
\begin{equation*}
6.4730\times10^{-2}.
\end{equation*}

\begin{figure}[H]
\centering\includegraphics[width=\linewidth]{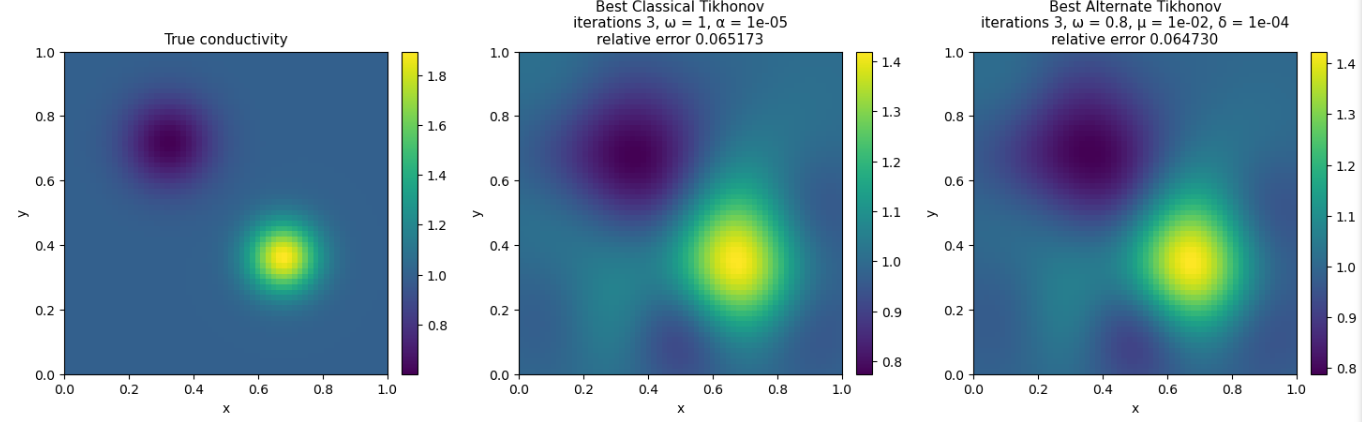}

\caption{Best reconstructions using each regularization method compared to true conductivity}
\end{figure}
The best classical Tikhonov result is obtained with
\begin{equation*}
n=3,\qquad\omega=1.0,\qquad\alpha=10^{-5},
\end{equation*}
which gives relative reconstruction error
\begin{equation*}
6.5173\times10^{-2}.
\end{equation*}
Which is slightly worse than the results from Alternative Tikhonov. Thus, Alternative Tikhonov regularization gives the best result in this experiment, although the difference between the best Alternative Tikhonov and best classical Tikhonov reconstructions is small.
A notable observation is that the best reconstructions occur after only a small number of outer iterations. So, increasing the number of iterations does not necessarily improve the result. This is due to overfitting of the regularization method, which begins to fit to error and noise leading to increasing error.
Future work could test these methods on more complicated conductivities, different noise levels, finer grids, and different choices of boundary inputs as well as tests on actual experimental data. It would also be useful to develop parameter choice rules for $\alpha$, $\mu$, $\delta$, $\omega$, and the stopping iteration $n$, rather than choosing them by a full parameter sweep.
Another important topic for further research is determining the best selections for boundary inputs $u_i|_{\partial\Omega}=f_i$ and $v_i|_{\partial\Omega}=g_i$. The number of inputs $N$ is worth investigating. It would also prove useful to put the boundary inputs in a form more accurate to experimental values. To this point performing actual experiments with an electrical impedence tomography device would be a crucial next step in this research. As well as investigating $f_i$ and $g_i$, in this paper we chose to set the adjoint equation $-\nabla\cdot(\gamma\nabla u)=p=0$, however it would be worthwhile to look into other options for $p$ close to 0.
Overall, the results show that Alternative Tikhonov regularization is a competitive alternative to classical Tikhonov regularization for the projected system $Ac=r$. It did not completely outperform classical Tikhonov in every setting, but it produced the lowest observed reconstruction error and gave a useful alternative way to balance residual size against coefficient size.

\clearpage
\section*{Addendum: Detailed Numerical Results}
\addcontentsline{toc}{section}{Addendum: Detailed Numerical Results}
\renewcommand{\arraystretch}{1.04}
\captionsetup[table]{font=footnotesize,skip=2pt}
\setlength{\intextsep}{6pt}
\newcommand{\resulthead}{\toprule Method & Iterations & Damping & Regularization parameters & Relative error\\\midrule}
\newcommand{\ct}{Classical Tikhonov}
\newcommand{\at}{Alternative Tikhonov}
\begin{table}[H]\centering\footnotesize
\textbf{Best Classical Tikhonov results by iterations}\par\vspace{0.5em}
\begin{tabular}{@{}lrrlr@{}}\resulthead
\ct & 3 & $\omega=1$ & $\alpha=1\times10^{-5}$ & 0.065173\\
\ct & 5 & $\omega=0.8$ & $\alpha=1\times10^{-5}$ & 0.065182\\
\ct & 4 & $\omega=1$ & $\alpha=1\times10^{-5}$ & 0.065183\\
\ct & 10 & $\omega=0.5$ & $\alpha=1\times10^{-5}$ & 0.065334\\
\ct & 20 & $\omega=0.5$ & $\alpha=3\times10^{-5}$ & 0.065598\\
\ct & 30 & $\omega=0.5$ & $\alpha=3\times10^{-5}$ & 0.065859\\
\bottomrule\end{tabular}
\caption{classical Tikhonov best results over n Iterations}
\end{table}

\begin{table}[H]\centering\footnotesize
\textbf{Best Classical Tikhonov results by damping}\par\vspace{0.5em}
\begin{tabular}{@{}lrrlr@{}}\resulthead
\ct & 3 & $\omega=1$ & $\alpha=1\times10^{-5}$ & 0.065173\\
\ct & 5 & $\omega=0.8$ & $\alpha=1\times10^{-5}$ & 0.065182\\
\ct & 10 & $\omega=0.5$ & $\alpha=1\times10^{-5}$ & 0.065334\\
\bottomrule\end{tabular}
\caption{classical Tikhonov best results over omega}
\end{table}

\begin{table}[H]\centering\footnotesize
\textbf{Best Classical Tikhonov results by alpha}\par\vspace{0.5em}
\begin{tabular}{@{}lrrlr@{}}\resulthead
\ct & 3 & $\omega=1$ & $\alpha=1\times10^{-5}$ & 0.065173\\
\ct & 10 & $\omega=1$ & $\alpha=3\times10^{-5}$ & 0.065585\\
\ct & 30 & $\omega=1$ & $\alpha=1\times10^{-4}$ & 0.065911\\
\ct & 30 & $\omega=1$ & $\alpha=3\times10^{-4}$ & 0.067343\\
\ct & 30 & $\omega=1$ & $\alpha=1\times10^{-3}$ & 0.069420\\
\ct & 30 & $\omega=1$ & $\alpha=3\times10^{-3}$ & 0.072521\\
\ct & 30 & $\omega=1$ & $\alpha=1\times10^{-2}$ & 0.076687\\
\bottomrule\end{tabular}
\caption{classical Tikhonov best results over alpha}
\end{table}

\begin{table}[H]\centering\footnotesize
\textbf{Best Alternative Tikhonov results by iterations}\par\vspace{0.5em}
\begin{tabular}{@{}lrrlr@{}}\resulthead
\at & 3 & $\omega=0.8$ & $\mu=1\times10^{-2}$, $\delta=1\times10^{-4}$ & 0.064730\\
\at & 4 & $\omega=0.8$ & $\mu=1\times10^{-2}$, $\delta=1\times10^{-3}$ & 0.064905\\
\at & 5 & $\omega=0.8$ & $\mu=1\times10^{-3}$, $\delta=3\times10^{-3}$ & 0.065187\\
\at & 20 & $\omega=0.5$ & $\mu=3\times10^{-2}$, $\delta=1\times10^{-3}$ & 0.065618\\
\at & 10 & $\omega=0.8$ & $\mu=3\times10^{-3}$, $\delta=3\times10^{-3}$ & 0.065622\\
\at & 30 & $\omega=0.8$ & $\mu=1\times10^{-1}$, $\delta=3\times10^{-4}$ & 0.065876\\
\bottomrule\end{tabular}
\caption{Alternative Tikhonov best results over n iterations}
\end{table}

\begin{table}[H]\centering\footnotesize
\textbf{Best Alternative Tikhonov results by damping}\par\vspace{0.5em}
\begin{tabular}{@{}lrrlr@{}}\resulthead
\at & 3 & $\omega=0.8$ & $\mu=1\times10^{-2}$, $\delta=1\times10^{-4}$ & 0.064730\\
\at & 3 & $\omega=1$ & $\mu=1\times10^{-2}$, $\delta=3\times10^{-4}$ & 0.064800\\
\at & 5 & $\omega=0.5$ & $\mu=1\times10^{-2}$, $\delta=3\times10^{-4}$ & 0.065193\\
\bottomrule\end{tabular}
\caption{Alternative Tikhonov best results omega}
\end{table}

\begin{table}[H]\centering\footnotesize
\textbf{Best Alternative Tikhonov results by mu}\par\vspace{0.5em}
\begin{tabular}{@{}lrrlr@{}}\resulthead
\at & 3 & $\omega=0.8$ & $\mu=1\times10^{-2}$, $\delta=1\times10^{-4}$ & 0.064730\\
\at & 5 & $\omega=0.8$ & $\mu=1\times10^{-3}$, $\delta=3\times10^{-3}$ & 0.065187\\
\at & 20 & $\omega=0.5$ & $\mu=3\times10^{-2}$, $\delta=1\times10^{-3}$ & 0.065618\\
\at & 10 & $\omega=0.8$ & $\mu=3\times10^{-3}$, $\delta=3\times10^{-3}$ & 0.065622\\
\at & 30 & $\omega=0.8$ & $\mu=1\times10^{-1}$, $\delta=3\times10^{-4}$ & 0.065876\\
\at & 30 & $\omega=1$ & $\mu=3\times10^{-1}$, $\delta=1\times10^{-4}$ & 0.071560\\
\bottomrule\end{tabular}
\caption{Alternative Tikhonov best results over mu}
\end{table}

\begin{table}[H]\centering\footnotesize
\textbf{Best Alternative Tikhonov results by delta}\par\vspace{0.5em}
\begin{tabular}{@{}lrrlr@{}}\resulthead
\at & 3 & $\omega=0.8$ & $\mu=1\times10^{-2}$, $\delta=1\times10^{-4}$ & 0.064730\\
\at & 3 & $\omega=1$ & $\mu=1\times10^{-2}$, $\delta=3\times10^{-4}$ & 0.064800\\
\at & 4 & $\omega=0.8$ & $\mu=1\times10^{-2}$, $\delta=1\times10^{-3}$ & 0.064905\\
\at & 5 & $\omega=0.8$ & $\mu=1\times10^{-3}$, $\delta=3\times10^{-3}$ & 0.065187\\
\at & 30 & $\omega=1$ & $\mu=1\times10^{-3}$, $\delta=1\times10^{-2}$ & 0.090462\\
\bottomrule\end{tabular}
\caption{Alternative Tikhonov best results over delta}
\end{table}

\begin{table}[H]\centering\footnotesize
\textbf{Best result by method}\par\vspace{0.5em}
\begin{tabular}{@{}lrrlr@{}}\resulthead
\at & 3 & $\omega=0.8$ & $\mu=1\times10^{-2}$, $\delta=1\times10^{-4}$ & 0.064730\\
\ct & 3 & $\omega=1$ & $\alpha=1\times10^{-5}$ & 0.065173\\
\bottomrule\end{tabular}
\caption{Best results for each method}
\end{table}

\begin{table}[H]\centering\footnotesize
\textbf{Top overall reconstruction results}\par\vspace{0.5em}
\begin{tabular}{@{}lrrlr@{}}\resulthead
\at & 3 & $\omega=0.8$ & $\mu=1\times10^{-2}$, $\delta=1\times10^{-4}$ & 0.064730\\
\at & 3 & $\omega=1$ & $\mu=1\times10^{-2}$, $\delta=3\times10^{-4}$ & 0.064800\\
\at & 4 & $\omega=0.8$ & $\mu=1\times10^{-2}$, $\delta=1\times10^{-3}$ & 0.064905\\
\at & 3 & $\omega=0.8$ & $\mu=1\times10^{-2}$, $\delta=3\times10^{-4}$ & 0.065014\\
\at & 3 & $\omega=1$ & $\mu=1\times10^{-2}$, $\delta=1\times10^{-4}$ & 0.065128\\
\ct & 3 & $\omega=1$ & $\alpha=1\times10^{-5}$ & 0.065173\\
\ct & 5 & $\omega=0.8$ & $\alpha=1\times10^{-5}$ & 0.065182\\
\ct & 4 & $\omega=1$ & $\alpha=1\times10^{-5}$ & 0.065183\\
\at & 5 & $\omega=0.8$ & $\mu=1\times10^{-3}$, $\delta=3\times10^{-3}$ & 0.065187\\
\at & 5 & $\omega=0.5$ & $\mu=1\times10^{-2}$, $\delta=3\times10^{-4}$ & 0.065193\\
\bottomrule\end{tabular}
\caption{Best results}
\end{table}
\end{document}